\documentclass[12pt]{amsart}

\usepackage{epsf,amssymb,times}

\newtheorem{thm}{Theorem}[section]
\newtheorem{lem}[thm]{Lemma}
\newtheorem{cor}[thm]{Corollary}

\newtheorem{rmk}[thm]{Remark}

\newtheorem{claim}[thm]{Claim}

\newcommand{\bdry}{\partial}

\newcommand{\s}{\vskip.1in}
\newcommand{\n}{\noindent}

\newcommand{\be}{\begin{enumerate}}
\newcommand{\ee}{\end{enumerate}}

\def\tb{\operatorname{tb}}

\def\Z{\hbox{$\mathbb Z$} }

\def\R{\hbox{$\mathbb R$} }

\def\dfn#1{{\em #1}}

\begin{document}

\title{Knots and Contact Geometry II: Connected Sums}
\date{May 28, 2002}

\author{John B.\ Etnyre}
\address{Department of Mathematics,
University of Pennsylvania,
209 South 33rd St.,
Philadelphia, PA 19104-6395}
\email{etnyre@math.upenn.edu}
\urladdr{http://math.upenn.edu/\char126 etnyre}

\author{Ko Honda}
\address{Department of Mathematics, University of Southern California, Los 
Angeles, CA 90089-1113} 
\email{khonda@math.usc.edu}
\urladdr{http://math.usc.edu/\char126 khonda}


\begin{abstract}
We study the behavior of Legendrian and transverse knots under the operation of 
connected sums. As a consequence we show that there exist Legendrian knots that
are not distinguished by any known invariant. Moreover, we classify 
Legendrian knots in some non-Legendrian simple knot types.
\end{abstract}

\maketitle


\section{Introduction}

The last few years have brought forth several advances in our 
understanding of Legendrian and transverse knots.  Roughly speaking, our 
knowledge has advanced on two fronts: via the holomorphic theory and via 
3-dimensional topology.  The most concrete realization of the holomorphic theory 
is the theory of Chekanov-Eliashberg contact homology invariants \cite{Chekanov, 
EGH}.  This theory yielded the first examples of nonisotopic Legendrian knots 
with the same classical invariants:  the topological type, the 
Thurston-Bennequin invariant, and the rotation number.  (There are also more 
computable variants derived from contact homology, 
such as the characteristic algebra of Ng \cite{Ng}.)  The 
purpose of these holomorphic invariants is to {\em distinguish}.  Their 
counterpart is 3-dimensional contact {\em topology}, which has the flavor of 
classical 3-dimensional cut-and-paste topology with a slight twist.  The main 
tool here is convex surface theory, introduced by Giroux \cite{Giroux}.  Using 
recent advances in convex surfaces, the authors completely classified Legendrian 
torus knots and Legendrian figure eight knots \cite{EH}.  A complete 
classification of Legendrian knots of a certain topological type implies the 
complete classification of transverse knots of the same topological type 
\cite{EH} (although not vice versa); hence transverse torus knots and transverse 
figure eight knots are classified (the predecessor to this result is 
\cite{Etnyre-transverse}).  More recently, Menasco \cite{Menasco} classified all 
transverse iterated torus knots by using the work of Birman-Wrinkle \cite{BW} 
which rephrased the classification question into a question in braid theory.

The goal of this paper is to prove a structure theorem for Legendrian 
knots, namely the behavior of Legendrian and transverse knots under 
the connected sum operation.  Our main theorem (Theorem~\ref{mainsum}) classifies 
Legendrian knots in a non-prime knot type, provided we understand the 
classification for the prime summands.  Theorem~\ref{mainsum}, in essence, is 
the relative version of Colin's gluing theorem for connected sums of tight contact 
manifolds \cite{Colin}.  One corollary of our main theorem is the existence of 
Legendrian knots which are not contact isotopic but are indistinguishable by all 
known invariants (including the holomorphic invariants).  Moreover, for any 
integer $m$, there exist Legendrian knots with identical invariants that are 
non-Legendrian-isotopic even after $m$ stabilizations.  Previously it was not 
known whether Legendrian knots (with identical invariants) became isotopic after 
one stabilization, largely due to the fact that the Chekanov-Eliashberg 
invariants vanish on stabilized Legendrian knots.  Theorem~\ref{mainsum} also 
implies the following: the connected sum of transversally simple knot types is 
transversally simple (see Section~\ref{main theorems} for definitions). 

The plan for the paper is as follows.  After reviewing some background 
(especially on connected sums of knots) in Section~\ref{background}, we give 
precise statements of the main theorem in Section~\ref{main 
theorems} and its applications in Section~\ref{appl}.  The main theorem is 
proved in Sections~\ref{section:standard} and \ref{Proofs}.

\section{Some background and notation}     \label{background}

We assume familiarity with basic notions in contact geometry, such as
characteristic foliations and convex surface theory.  This can be found in 
\cite{Etnyre-GA} (see also \cite{Aebischer, Eliashberg92, Giroux}).  As this 
paper is a sequel to \cite{EH} we assume the reader is familiar with its 
contents.  In particular,  Sections~2 and 3 of \cite{EH} are foundational 
and develop the necessary terminology and background on Legendrian knots and 
transversal knots.

In this paper, our ambient 3-manifolds and knots are {\em oriented}, and ``knot 
types'' are oriented knot types.  Let $M_1$ and $M_2$ be 3-manifolds.  We first 
describe the {\em connected sum} of two (topological) knots $K_1\subset M_1$ and 
$K_2\subset M_2$.  Let $B_i$ be an open ball in $M_i$ that intersects $K_i$ in 
an unknotted arc $\alpha_i.$  Let $f:\partial (M_1\setminus B_1)\to \partial 
(M_2\setminus B_2)$ be an orientation-reversing diffeomorphism which sends 
$K_1\cap \partial (M_1\setminus B_1)$ to $K_2\cap \bdry (M_2\setminus B_2).$  
(Here, $X\setminus Y$ denotes the metric closure of the complement of $Y$ in 
$X$.)  Now the connected sum of $M_1$ and $M_2$ is \[M_1\# M_2=(M_1\setminus 
B_1)\cup_f (M_2\setminus B_2)\] and the {\dfn{connected sum} of $K_1$ and $K_2$ 
in $M_1\# M_2$ is \[K_1\# K_2 = (K_1\setminus \alpha_1)\cup (K_2\setminus 
\alpha_2).\] Note that there are two possible identifications of $K_1\cap 
\partial (M_1\setminus B_1)$ with $K_2\cap \bdry (M_2\setminus B_2)$ --- we 
choose the one which induces a coherent orientation on $K_1\# K_2$.  It is an 
easy exercise to see that $K_1\# K_2$ is well-defined and its topological 
type is independent of the choices of $B_i$ and $f$.    

If $K_1, K_2\subset S^3$, we can interpret the connected sum operation as 
happening entirely in $S^3$, since $S^3\# S^3=S^3$.  In particular, fix a  
2-sphere $S$ in $S^3$ that splits $S^3$ into two balls $B_1$ and $B_2.$  Then 
isotop $K_1$ so that it intersects $B_2$ in an unknotted arc and isotop $K_2$ 
so that it intersects $B_1$ in an unknotted arc.  Moreover, we can arrange for 
$K_1$ and $K_2$ to intersect $S$ at the same points.  We then define $K_1\# K_2= 
(K_1\cap B_1)\cup (K_2\cap B_2).$   This clearly is an equivalent definition of 
the connected sum in $S^3.$  From this definition it is easy to arrive at the 
familiar diagrammatic picture of a connected sum. See Figure~\ref{fig:sump}.

\begin{figure}[ht]
	{\epsfxsize=3in \centerline{\epsfbox{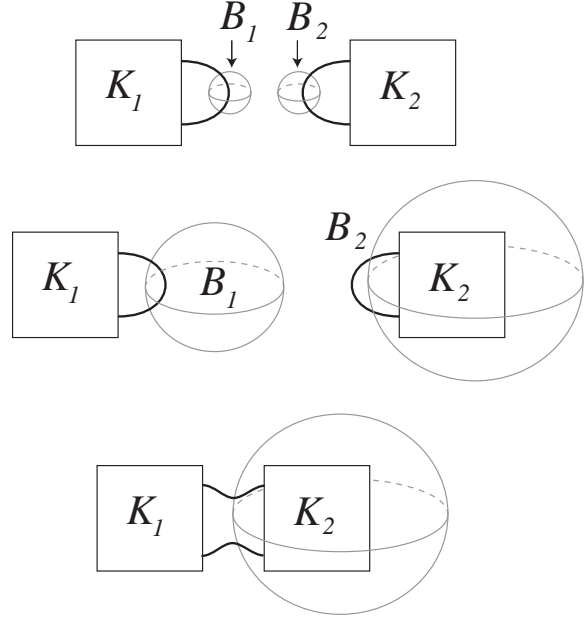}}
	\caption{Diagrammatic connected sums.}
	\label{fig:sump}}
\end{figure}

A knot $K$ in $S^3$ is \dfn{prime} if $K=K_1\# K_2$ implies that either $K_1$ or 
$K_2$ is the unknot.  Any knot $K\subset S^3$ admits a unique (minimal)
decomposition into prime pieces, {\em i.e.}, $K=K_1\# \dots \# K_n$ with $K_i$, 
$i=1,\dots,n$, prime.  This decomposition can be achieved by finding a 
collection $\{S_1,\dots, S_{n-1}\}$ of disjoint 2-spheres in $S^3$ that each 
intersects $K$ in two points.  Given such a {\em separating sphere} $S_i$, we 
may reverse the procedure described in the preceding paragraph to write the knot 
as the connected sum of two other knots.

Although the collection $\{K_1,\dots, K_n\}$ is unique up to isotopy, the 
collection $\{S_1,\dots, S_{n-1}\}$ of separating spheres is not.  To avoid 
confusion in what follows, whenever we decompose a knot in $S^3$, we will be 
doing so with respect to a {\em fixed} collection of separating spheres 
$\{S_1,\dots, S_{n-1}\}$.  Moreover, we take $K=K_1\#\dots\#K_n$ to mean the 
following:  using the same notation from the second paragraph of this section,  
we glue the $K_i\setminus \alpha_i$ together so that the endpoint of 
$K_i\setminus \alpha_i$ connects to the initial point of $K_{i+1}\setminus 
\alpha_{i+1}$ modulo $n$.  (This makes sense since the $K_i$ are oriented.) This 
way, the $K_1,\dots,K_n$ are cyclically strung together in order.

Let $\mathcal{K}$ be a topological knot type in a 3-manifold $M$, {\em i.e.}, an 
equivalence class of (topologically) isotopic knots.  Define 
$\mathcal{L}(\mathcal{K}, M, \xi)$ to be the set of isotopy classes of 
Legendrian knots in $(M,\xi)$ of type $\mathcal{K}.$  If the contact manifold 
$(M, \xi)$ is implicit, then we write $\mathcal{L}(\mathcal{K})$ instead of 
$\mathcal{L}(\mathcal{K}, M, \xi).$

\section{The Main Theorem} \label{main theorems}

We first explain Colin's gluing theorem~\cite{Colin}.  Denote the space of tight 
contact 2-plane fields on a 3-manifold $M$ by $Tight(M)$.  Then we have the 
following:
                                          
\begin{thm}[Colin] \label{thm:colin}
Given two 3-manifolds $M_1, M_2$, there is an isomorphism
$$\pi_0(Tight(M_1))\times \pi_0(Tight(M_2))\stackrel 
\sim\longrightarrow\pi_0(Tight(M_1\# M_2)).$$ 
\end{thm}

Let $(M_i, \xi_i)$, $i=1,2$, be two tight contact manifolds.  Choose $p_i\in 
M_i$ as well as a standard contact 3-ball $B_i$ with coordinates $(x,y,z)$ about 
$p_i$ so that the contact structure is given by $dz+xdy=0$.  After possibly 
perturbing the boundary of $B_i$, there is an orientation-reversing map $f$ from 
$S_1=\partial (M_1\setminus B_1)$ to $S_2=\partial (M_2\setminus B_2)$ that 
takes the characteristic foliation of $S_1$ to that of $S_2$.  According to 
Colin's theorem, the contact structure $\xi$ induced on \[M=M_1\# M_2= 
(M_1\setminus B_1) \cup_f (M_2\setminus B_2)\] is tight, and is independent of 
the choice of $B_i$, $p_i$, and $f$, up to isotopy.  Moreover, every tight $\xi$ 
on $M$ arises, up to isotopy, from a unique pair $(\xi_1,\xi_2)$ of tight 
contact structures.

Let us now explain the Legendrian connected sum operation, which is a relativized 
version of Colin's connected sum operation.  In each $(M_i,\xi_i)$, choose an 
oriented Legendrian knot $L_i$ and a point $p_i\in L_i.$   Normalize the 
standard contact 3-ball $B_i$ so that $L=B_i\cap \text{$y$-axis}$, and 
further require that $f$ maps $L_1\cap S_1$ to $L_2\cap S_2$ as oriented 
manifolds. Then we obtain the Legendrian knot $L= L_1\# L_2 \subset M$, which 
is called the \dfn{connected sum} $L_1\#L_2$ of $L_1$ and $L_2.$ 

\begin{lem}
The connected sum of two Legendrian knots does not depend on the points $p_i$, 
the balls $B_i$, or $f$ used in the definition. 
\end{lem}

Although this lemma is not difficult to prove, we defer the proof until 
Section~\ref{Proofs}.  See \cite{Chekanov} for a diagrammatic proof.

Given a nullhomologous Legendrian knot $L$, we can define its Thurston-Bennequin 
invariant $\tb(L)$ and it rotation number $r(L)$.  (For more details, consult
\cite{Etnyre-GA} and \cite{EH}, for example.)  We then have:

\begin{lem}
If $L_1$ and $L_2$ are two nullhomologous Legendrian knots, then
\begin{equation}
	\tb(L_1\# L_2) = \tb (L_1) + \tb (L_2) + 1, 
\end{equation}
and
\begin{equation}
	r(L_1\# L_2) = r (L_1) + r (L_2). 
\end{equation}
\end{lem}

This lemma easily follows from the facts that $\tb$ and $r$ can be computed from 
the characteristic foliation of a Seifert surface of a knot (see \cite{EF}) and 
that we can control the characteristic foliation on the Seifert surface for 
$L_1\# L_2$ in terms of the foliations on the surfaces for $L_1$ and $L_2.$

We denote by $S_\pm(L)$ the $\pm$ stabilization of the Legendrian knot $L$.  
Recall that this is a procedure to reduce $\tb$ of a Legendrian knot by $1$ (see 
\cite{EH}) and diagrammatically corresponds to ``adding kinks'' to $L.$   The 
following is our main theorem.

\begin{thm}\label{mainsum}
Let $\mathcal{K}= \mathcal{K}_1\#\dots \#\mathcal{K}_n$ be a 
connected sum decomposition of a topological knot type $\mathcal{K}\subset 
(M,\xi)$ into prime pieces $\mathcal{K}_i\subset (M_i,\xi_i)$, where $(M,\xi)= 
(M_1,\xi_1)\#\dots\#(M_n,\xi_n)$ is tight.  Then the map 	\begin{equation} 
C:\left({\mathcal{L}(\mathcal{K}_1)\times\dots\times \mathcal{L}(\mathcal{K}_n) 
\over \sim} \right) \longrightarrow 
\mathcal{L}(\mathcal{K}_1\#\dots \#\mathcal{K}_n) 	
\end{equation} 
given by $(L_1,\dots,L_n)\mapsto L_1\# \dots \# L_n$ is a bijection. 
Here the equivalence relation $\sim$ is of two types:
\be
\item $(L_1,\dots,S_{\pm}(L_i),L_{i+1},\dots,L_n)\sim
(L_1,\dots,L_i,S_{\pm}(L_{i+1}),\dots,L_n),$
\item $(L_1,\dots,L_n)\sim \sigma(L_1,\dots,L_n)$, where $\sigma$ is a
permutation of the $\mathcal{K}_i\subset (M_i,\xi_i)$ such that 
$\sigma(M_i,\xi_i)$ is isotopic to $(M_i,\xi_i)$ and 
$\sigma(\mathcal{K}_i)=\mathcal{K}_i$. 
\ee
\end{thm}

Theorem~\ref{mainsum} will be proved in Section~\ref{Proofs}.   We now 
discuss its consequences.  Let $\overline{\tb}(\mathcal{K})$ denote the maximal 
Thurston-Bennequin invariant over elements in $\mathcal{L}(\mathcal{K})$.  Then,

\begin{cor}
$\overline{\tb} (\mathcal{K}_1\# \mathcal{K}_2) = \overline{\tb} (\mathcal{K}_1) 
+ \overline{\tb} (\mathcal{K}_2) + 1.$
\end{cor}

This corollary was independently proven by Torisu in \cite{Torisu}.

Theorem~\ref{mainsum} takes on a particularly simple form when one restricts 
attention to maximal Thurston-Bennequin knots.

\begin{cor}   \label{maximal}
If $L$ is a Legendrian knot which is a maximal $\tb$ representative of  
$\mathcal{L}(\mathcal{K})$, then $L$ admits a unique prime decomposition
(modulo potential permutations). 
\end{cor}

Recall the strategy described in \cite{EH} for classifying Legendrian knots of a 
given a topological knot type.  The idea was to (1) prove that Legendrian knots 
in a particular knot type would always destabilize to a knot with maximal $\tb$ 
and then (2) classify all Legendrian knots of this type with maximal $\tb$.  
When executing the final stage of this strategy, Corollary~\ref{maximal} is 
useful because it allows us to concentrate on prime knots.

\section{Applications}   \label{appl}

In discussing the applications of Theorem~\ref{mainsum}, we restrict our 
attention to Legendrian knots in the standard tight contact $(S^3,\xi_{S^3})$ or 
$(\R^3, \xi_{\R^3})$.   Similar results hold in other manifolds.

First we reformulate the connected sum operation.  Let $\mathcal{K}_1$ and 
$\mathcal{K}_2$ be two topological knot types.  Given 
$L_i\in\mathcal{L}(\mathcal{K}_i)$, $i=1,2$, we define their \dfn{connected sum} 
as follows: let $B_i$ be a standard contact 3-ball about a point $p_i$ on 
$L_i.$  By Eliashberg's classification theorem of tight contact structures on 
the 3-ball \cite{Eliashberg92}, there is a contact isotopy $\phi_t:S^3\to S^3, 
t\in[0,1],$  from $\phi_0=\hbox{id}_{S^3}$ to $\phi_1$ which takes $B_1$ to 
$S^3\setminus B_2.$  Moreover, it is easy to arrange $L_2$ and 
$\phi_1(L_1)$ to intersect $\partial B_2$ in the same two points.  We may now 
define $L_1 \# L_2$ to be the Legendrian knot $(L_2\cap(S^3\setminus B_2)) \cup 
(\phi_1(L_1)\cap B_2)$.  We leave it as a simple exercise to check that this 
definition of the connected sum of knots is equivalent to the one given above. 
It has the advantage of being done ambiently, {\em i.e.}, we do not take 
connected sums of the ambient manifolds, only of the knots.

Since any Legendrian isotopy can be assumed to miss a preassigned point, the 
classification of Legendrian knots in $(S^3,\xi_{S^3})$ and in 
$(\R^3,\xi_{\R^3})$ are equivalent.  Moreover, there is a convenient 
diagrammatic description of Legendrian knots in $\R^3$ in terms of front 
projections (for example, see \cite{EH}). 
\begin{figure}[ht]
	{\epsfxsize=5in \centerline{\epsfbox{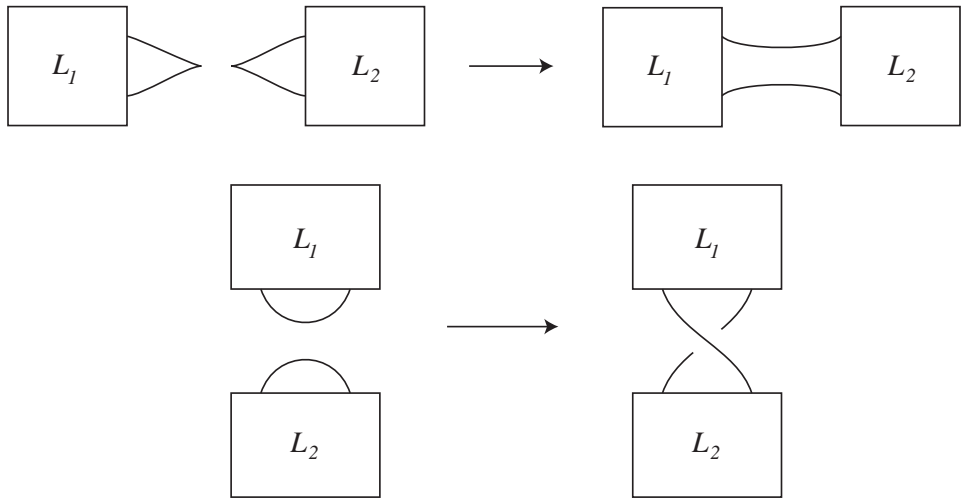}}
	\caption{Two diagrammatic versions of Legendrian connected sum.} 	
	\label{fig:diag}}
\end{figure}
Figure~\ref{fig:diag} indicates two ways in which the ambient connected sum 
(described in the previous paragraph) can be done in terms of the front 
projections of Legendrian knots.

Perhaps the most interesting application of Theorem~\ref{mainsum} is towards the 
construction of topological knot types which are not Legendrian simple.  Recall 
that a topological knot type $\mathcal{K}$ is said to be {\em Legendrian 
simple} if Legendrian knots in $\mathcal{K}$ are classified by the 
Thurston-Bennequin invariant and the rotation number.  The first 
non-Legendrian-simple knot type was discovered in \cite{Chekanov} and, since 
then, many similar examples have been found.  All examples to date have used 
contact homology (in one form or another) to distinguish the Legendrian knots.  
Although contact homology provides an intriguing way of distinguishing 
Legendrian knots, it currently does not provide much geometric insight into why 
Legendrian knots are different.

\begin{thm}\label{thm:examples} 	
Given two positive integers $m$ and $n$, there is a knot type $\mathcal{K}$ and 
distinct Legendrian knots $L_1,\dots, L_n$ in $\mathcal{L}(\mathcal{K})$ which 
have the same Thurston-Bennequin invariant and rotation number, and remain 
distinct even after $m$ stabilizations (of any type). 
\end{thm}

Theorem~\ref{thm:examples} follows from Theorem~\ref{mainsum} and the 
classification of Legendrian torus knots (Theorem~\ref{negtorus} below) from  
\cite{EH}.  Recall that $\mathcal{K}_{p,q}$ is a {\em $(p,q)$-torus knot} if any 
element of $\mathcal{K}_{p,q}$ can be isotoped to sit on a {\em standardly 
embedded} torus $T$ in $S^3$ as a $(p,q)$-curve.  Here we say a torus $T\subset 
S^3$ is {\em standardly embedded} if it is oriented and $S^3\setminus T= N_1\cup 
N_2$, where $N_i$, $i=1,2$, are solid tori with $\bdry N_1=T$ and $\bdry 
N_2=-T$.  Now, there exists an oriented identification $T\simeq \R^2/\Z^2$ where 
the meridian of $N_1$ corresponds to $\pm(1,0)$ and the meridian of $N_2$ to 
$\pm(0,1)$.

\begin{thm}\label{negtorus}
Legendrian knots in $\mathcal{L}(\mathcal{K}_{p,q})$ are determined by their knot type, 
Thurston-Bennequin invariant and rotation number.  If $p<0$ and $-p>q>0,$ then 
$\overline{\tb}(\mathcal{K}_{p,q})=pq$ and the corresponding values of $r$ are  
$$r(K)\in\{ \pm(\vert p\vert-\vert q\vert - 2qk) : k\in\mathbb{Z},
0\leq k<  \frac{\vert p\vert-\vert q\vert}{\vert q\vert}\}.$$
Moreover, all other Legendrian knots in this knot type are obtained by 
stabilization.
\end{thm}

If we plot the possible values of $\tb$ and $r$ for a negative torus knot, we 
obtain a picture similar to that of Figure~\ref{mountain}. 

\begin{figure}[ht]
	{\epsfxsize=5in \centerline{\epsfbox{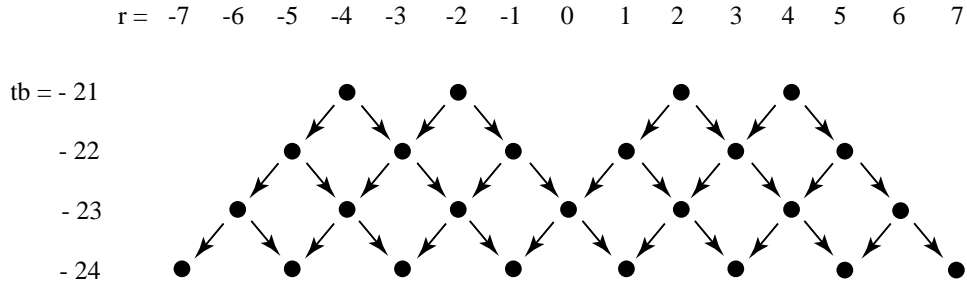}}
	\caption{Some possible $\tb$ and $r$'s for the $(-7,3)$ torus knot.}
	\label{mountain}}
\end{figure}

\begin{proof}[Proof of Theorem~\ref{thm:examples}.]
Let $p=-(4n+1)s -1$ and $q=2s$, with $s$ an even number greater than $m+1$.  
Then, according to Theorem~\ref{negtorus}, there are $4n$ Legendrian knots in 
$\mathcal{L}(\mathcal{K}_{p,q})$ with maximal $\tb=pq$ and distinct 
rotation numbers $-(4n-3)s+1,\dots, (4n-5)s+1,(4n-1)s+1$ and $-(4n-1)s-1, 
-(4n-5)s-1,\dots, (4n-3)s -1$.  Let $L_r\in \mathcal{L}(\mathcal{K}_{p,q})$ with 
$\tb=pq$ and rotation number $r$.  For $k=0,\dots, 2n-1,$ let 
$L^k=L_{(4(n-k)-1)s+1}\# L_{-(4(n-k)-1)s-1}.$  Note that all the $L^k$ are 
topologically isotopic, have the same $\tb=2pq+1$ and $r=0$, yet are not 
Legendrian isotopic by Theorem~\ref{mainsum}.  See Figure~\ref{fig:ex}.  Since 
the spacings in $r$ between adjacent maximal $\tb$ representatives are at least 
$2m$ by our choice of $p$ and $q$, the $L^k$ remain distinct even after $m$ 
stabilizations. \end{proof}

\begin{figure}[ht]
	{\epsfxsize=4.5in \centerline{\epsfbox{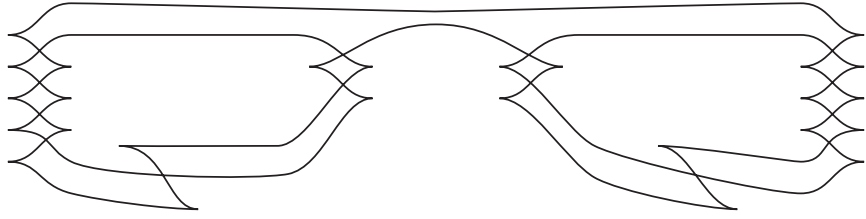}}
	\caption{When $m=0$ in Theorem~\ref{thm:examples}, one can use
		$(p,q)=(-(2n+1),2)$-torus knots. Here is one of those knots when $p=-7.$}
	\label{fig:ex}}
\end{figure}

\begin{rmk}
{\em The Legendrian knots $L^k$ appearing in the proof of Theorem~\ref{thm:examples} 
remain distinct after $m$ stabilizations.  However, it is well-known that the 
Chekanov-Eliashberg contact homology invariants are unable to distinguish 
stabilized knots (because the invariants vanish).  Thus these are the first 
examples of Legendrian knots which are not distinguished by the known 
holomorphic invariants.  We also note that the examples in 
Theorem~\ref{thm:examples} have nontrivial contact homology.  We are unable to 
determine if these invariants are the same or not, but all easily computable 
invariants derived from contact homology are the same for these examples. }
\end{rmk}

\begin{rmk}{\em 
Using Theorems~\ref{mainsum} and \ref{negtorus}, we can classify Legendrian 
knots isotopic to (multiple) connected sums of torus knots.  This is the first 
classification of Legendrian knots in a non-Legendrian-simple knot type. }
\end{rmk}

\begin{rmk}\label{notsimplerem}{\em 
Observe that the connected sums of torus knots are fibered knots.  Thus we have 
examples of non-Legendrian-simple fibered knots, contrary to a (perhaps overly 
optimistic) conjecture that fibered knots are Legendrian simple.}
\end{rmk}

We end by observing that, while the connected sum of Legendrian simple knot 
types need not be Legendrian simple, the connected sum of transversally simple 
knot types is always transversally simple.  Here a knot type $\mathcal{K}$ is 
{\em transversally simple} if transversal knots in $\mathcal{K}$ are determined 
by their self-linking number.

\begin{thm}
The connected sum of transversally simple knot types is transversally simple. 
\end{thm}

\begin{proof}
Recall that, according to Theorem~2.10 of \cite{EH}, a knot type is 
transversally simple if and only if it is stably simple.  A knot type 
$\mathcal{K}$ is {\em stably simple} if any two knots in 
$\mathcal{L}(\mathcal{K})$ for which $s=\tb-r$ agree are Legendrian isotopic 
after some number of negative stabilizations.

Now assume that $\mathcal{K}_1$ and $\mathcal{K}_2$ are stably simple knot 
types.  Let $L_1,L'_1\in \mathcal{L}(\mathcal{K}_1)$ and $L_2, L_2'\in 
\mathcal{L}(\mathcal{K}_2)$ such that $s(L_1\# L_2)=s(L_1'\# L_2').$  It follows 
that $s(L_1)= s(L_1')+2n$ and $s(L_2)=s(L_2')-2n$ for some integer $n$, which 
we may take to be $\geq 0$.  Since $\mathcal{K}_1$  and $\mathcal{K}_2$ are 
stably simple, there exist $m_1$ and $m_2$ such that $S^{m_1}_-\circ S_+^n(L_1)$ 
is Legendrian isotopic to $S^{m_1}_-(L'_1)$ and $S^{m_2}_-\circ S^n_+ (L_2')$ is 
Legendrian isotopic to $S_-^{m_2}(L_2).$  Thus 
\begin{align*} S^{m_1+m_2}_-(L_1\# L_2)&=(S_-^{m_1}(L_1))\# 
(S_-^{m_2}(L_2))= (S_-^{m_1}(L_1))\# (S^{m_2}_-\circ S_+^n (L_2'))\\ &= 
(S_-^{m_1}\circ S^n_+ (L_1))\# (S^{m_2}_-(L_2'))= (S_-^{m_1}(L_1'))\# 
(S^{m_2}_-(L_2'))\\ &= S^{m_1+m_2}_-(L_1'\# L_2').
\end{align*}
This proves that $\mathcal{K}_1\# \mathcal{K}_2$ is stably simple.
\end{proof}

\begin{rmk}{\em In contrast to the situation for Legendrain knots discussed in 
Remark~\ref{notsimplerem}, it still does not seems unreasonable
to believe that fibered knots are transversely simple. See also \cite{BW, Menasco}.
}\end{rmk}

\section{The main technical result}\label{section:standard}

Given a Legendrian knot $L$ in a tight contact manifold $(M,\xi)$, we may always 
find a sufficiently small tubular  neighborhood $N$ of $L$ such that $T=\partial 
N$ is a convex torus with dividing set $\Gamma_T$ consisting of two parallel, 
homotopically nontrivial dividing curves.  We make an oriented identification 
$T\simeq \R^2/\Z^2$ with coordinates $(\mu,\lambda)$, so that the 
$\mu$-direction is the meridional direction and $\lambda$-direction is the 
longitudinal direction given by a Seifert surface.  (Note that this convention 
is different from the usual Dehn surgery convention.) The slope of a 
homotopically nontrivial closed curve on $T$ will be given in the 
$\mu\lambda$-coordinates.  With respect to these coordinates, the slope of 
$\Gamma_T$ is $\frac{1}{\tb(L)}.$  Using the Legendrian Realization Principle we 
may arrange, and shall always assume, that $T$ is in {\em standard form} and the 
ruling slope on $T$ is $0.$ 

An embedded sphere $S$ in $M$ that intersects $L$ transversely in exactly two 
points and separates $M$ will be called a \dfn{separating sphere} for $(M,L).$
Given a separating sphere $S$, let $M\setminus S= M^o_1 \sqcup M^o_2$ and 
$L^o_i=(L\setminus S)\cap M^o_i$, $i=1,2$.    We call $S$ a \dfn{trivial} 
separating sphere if one of the $M_i^o$ is a 3-ball and $L^o_i$ is an unknotted 
arc in $M^o_i=B^3.$  The separating sphere $S$ can be isotoped so that $S\cap T$ 
($T=\bdry N$) consists of two ruling curves.  We may further isotop $S$, 
relative to $S\cap T$, so that $S$ becomes convex and the annular component of 
$S\setminus T$ admits a ruling by closed Legendrian curves parallel to the 
boundary of the annulus.  Such an $S$ will be called a {\em standard convex 
separating sphere}.

We now introduce a standard object to cap off our cut-open manifold/knot pairs 
$(M_i^o, L_i^o).$  To this end, let $N_s$ be a convex tubular neighborhood of 
the $y$-axis in the standard tight contact $(\R^3,\xi_{std})$ given by the 
1-form $dz+xdy$.  We can assume the dividing curves on $N_s$ consist of two 
lines parallel to the $y$-axis and arrange the ruling curves to be all 
meridional.  Now let $B_s$ be a 3-ball about the origin with convex $\bdry B_s$, 
such that $\partial B_s\cap \partial N_s$ consists of two ruling curves. 
Finally, let $L_s= \text{$y$-axis}\cap B_s.$  We call the pair 
$((B_s,\xi_{std}|_{B_s}), L_s)$, consisting of the tight contact manifold 
$(B_s,\xi_{std}|_{B_s})$ and the Legendrian arc $L_s$, the {\em standard trivial 
pair}.

Given a convex separating sphere $S$ as above, we can apply the Giroux 
Flexibility Theorem so that the characteristic foliations on $S$ and $\partial 
B_s$ agree.  For each $i=1,2$, we then glue $(B_s, L_s)$ to $(M_i^o, L_i^o)$ to 
get a closed contact manifold $(M_i,\xi_i)$ and a Legendrian knot $L_i\subset 
M_i.$  The following is a consequence of Theorem~\ref{thm:colin}.

\begin{cor}\label{mainambient}  
$(M_i,\xi_i)$ is tight, and, up to isotopy, does not depend on the choice of 
convex separating sphere $S$ (provided the topological type is preserved) or on 
the gluing map.  
\end{cor}

We now consider the relative version of the corollary which takes into account 
the splitting of the Legendrian knot $L\subset M$.  We then have:

\begin{thm}\label{maintech}
Let $((M,\xi),L)$ be a tight contact manifold together with a Legendrian knot 
$L\subset M$, and let $S$, $S'$ be (smoothly) isotopic standard convex 
separating spheres.  Let $((M_i,\xi_i),L_i)$ (resp.\ $((M_i,\xi_i),L_i')$), 
$i=1,2$, be the glued-up manifolds together with Legendrian knots, obtained by 
cutting $M$ along $S$ (resp.\ $S'$) and gluing in copies of the standard trivial 
pair.  Then there exists a sequence $(L_1^0,L_2^0)=(L_1,L_2), 
(L_1^1,L_2^1),\dots, (L_1^k, L_2^k)=(L_1',L_2')$, where:
\be
\item $L_i^j$ is a Legendrian knot in $(M_i,\xi_i)$ isotopic to, but not 
necessarily contact isotopic to, $L_i$, $i=1,2$, 
\item $(L_1^{j+1},L_2^{j+1})$ is obtained from $(L_1^j,L_2^j)$ by performing one 
of the following: 
\be 
\item[(i)] Legendrian isotopy,
\item[(ii)] $L_1^{j+1}= S_\pm(L_1^j)$ and $L_2^{j+1}=(S_\pm)^{-1}(L_2^j)$, or
\item[(iii)] $L_1^{j+1}= (S_\pm)^{-1}(L_1^j)$ and $L_2^{j+1}=S_\pm(L_2^j)$. 
\ee
Here, $(S_\pm)^{-1}$ indicates destabilization.
\ee
\end{thm}

In other words, the Legendrian knots $L_i$ and $L_i'$ which arise from isotopic 
but not contact isotopic separating spheres differ only by successively shifting 
Legendrian stabilizations from one side to the other.

\s\n
The remainder of this section is devoted to the proof of Theorem~\ref{maintech}.  
The proof is essentially a concrete application of the state traversal 
technique.

\s\n
{\bf Step 1.}  Let $\xi$ be a $[0,1]$-invariant tight contact structure on 
$A=S^2\times[0,1]$, viewed as a neighborhood of $\partial B_s$ sitting in 
$(\R^3,\xi_{std}).$  It follows from the proof of Eliashberg's classification of 
tight contact structures on the 3-ball \cite{Eliashberg92} that $\xi$ is the 
unique (up to isotopy rel boundary) tight contact structure on $A$, with the 
given characteristic foliation on the boundary. 

The intersection of $A$ with the $y$-axis has two components $L_+$ and $L_-$; 
these are Legendrian arcs running between the two boundary components of $A.$  
Let $L^{m,n}_\pm=S_+^m\circ S_-^n(L_\pm).$

\begin{figure}[ht]
	{\epsfxsize=2in \centerline{\epsfbox{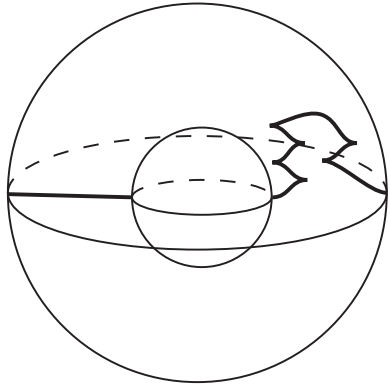}}
	\caption{The arcs $(L_-, L^{2,1}_+)$ in $A.$}
	\label{ex}}
\end{figure}

\begin{lem}\label{facts}
Let $L'_+$ and $L_-'$ be Legendrian arcs in $A$ which have the same endpoints 
as $L_+$ and $L_-$, respectively, and such that $L'_+\sqcup L'_-$ is (smoothly) 
isotopic to $L_+\sqcup L_-$ inside $A$, rel $\bdry A$.  Then, after applying a 
contactomorphism which is isotopic to the identity through an isotopy which 
fixes only one of the boundary components, $L'_-\sqcup L_+'$ is Legendrian 
isotopic to $L_-\sqcup L_+^{m,n}$ for some uniquely determined $m$ and $n.$ 
\end{lem}

\begin{proof}[Proof of Lemma~\ref{facts}.]
We define the {\em twisting number} (or the {\em relative Thurston-Bennequin 
invariant}) $tw(L'_\pm)$ to be the difference between the contact framings of 
$L'_\pm$ and $L_\pm$.  Note that the well-definition of $tw(L'_\pm)$ follows 
from the fact that $L'_\pm$ and $L_\pm$ have the same endpoints.  

Now, let $g$ be the diffeomorphism of $A=S^2\times[0,1]$ which rotates the 
sphere $S^2\times\{t\}$ around the axis provided by $L_+\sqcup L_-$ ({\em i.e.}, 
the $y$-axis) by $2\pi k t$, where $k$ is chosen so that $tw(g(L'_-))=0$ with 
respect to $g_*\xi$.  Here $g_*\xi$ is isotopic to $\xi$ rel boundary by 
the uniqueness of tight contact structures on $A$ with fixed boundary.

Observe that $tw(L'_+) \leq 0$, since otherwise we could use $L'_-$, $L'_+$ and 
arcs on $\partial A$ to construct a Legendrian unknot in $A$ with nonnegative 
Thurston-Bennequin invariant, which would contradict tightness.  Let $m$ and $n$ 
be nonnegative integers which satisfy $m+n=-tw(L'_+).$  (The precise values of 
$m$ and $n$ are to be determined later.)  Hence $tw(L_+^{m,n})=tw(L'_+)$.  
Next, there exists an isotopy $f$ of $A$ rel $\bdry A$ which takes $L'_-\sqcup 
L'_+$ to $L_-\sqcup L_+^{m,n}$.  Since Legendrian curves and their standard 
tubular neighborhoods are interchangeable for all practical purposes, we may 
assume that $f$ is a contactomorphism from the neighborhood $U$ of $L'_-\sqcup 
L'_+$ onto its image. 

It remains to extend $f$ to a contactomorphism on all of $A$ or, equivalently, 
match up two tight contact structures on the solid torus $A\setminus U$.  For 
this, we apply the classification of tight contact structures on solid tori 
\cite{Giroux2,Honda}.  The boundary slope for both tight contact structures on 
the solid torus $A\setminus U$ is $-(m+n)-1$.  By the classification, there 
exists a bijection between nonnegative integer pairs $(m,n)$ with $m+n= 
-tw(L'_+)$ and tight contact structures on $A\setminus U$ with slope 
$tw(L'_+)-1$ and two dividing curves on the boundary.  Hence there is a unique 
choice of $m$, $n$ so that the two tight contact structures on $A\setminus U$ 
are contact isotopic rel boundary.  \end{proof}

\s\n
{\bf Step 2.} 
We now establish Theorem~\ref{maintech} under an extra hypothesis on the spheres 
$S$ and $S'.$ 

\begin{claim}\label{lem2}
Theorem~\ref{maintech} holds if $S$ and $S'$ are disjoint and cobound a region 
diffeomorphic to $S^2\times[0,1].$ 
\end{claim}

\begin{proof}[Proof of Claim~\ref{lem2}.]
Let $A'\subset M$ be the region between $S$ and $S'$, and let $M^c_1$ and 
$M^c_2$ be components of $M\setminus A'$ so that: 
\begin{align*}
M_1&=M_1^c\cup B_s,\\
M_2&=M_2^c\cup A' \cup B_s,\\
M'_1&= M_1^c\cup A' \cup B_s, \text{ and}\\
M'_2&= M_2^c \cup B_s,
\end{align*}
where $B_s$ is the standard contact 3-ball.

Let $L_i^c=L\cap M_c^i$, $i=1,2$ and $L'=L \cap A'.$  Thus the Legendrian arcs 
under consideration are:
\begin{align*}
L_1&=L_1^c\cup L_s,\\
L_2&=L_2^c\cup L' \cup L_s,\\
L'_1&=L_1^c\cup L'\cup L_s, \text{ and}\\
L'_2&=L_2^c\cup L_s,
\end{align*}
where $L_s$ is the standard Legendrian arc in $B_s.$

Observe that $B_s$ is contactomorphic to $B_s\cup A$, where $A=S^2\times [0,1]$ 
with the $[0,1]$-invariant tight contact structure.  Therefore we may think of 
$M_1$ and $M'_2$ as composed of the appropriate $M_i^c, B_s$ and $A.$  Let 
$f:M_1\stackrel\sim\to M'_1$ be the diffeomorphism which sends $M_1^c\subset 
M_1$ to $M_1^c\subset M'_1$  by the identity, $B_s\subset M_1$ to $B_s\subset 
M'_1$ by the identity, and $A$ to $A'$ by a diffeomorphism preserving the 
characteristic foliation on the boundary.  It is easy to arrange for the 
diffeomorphism from $A$ to $A'$ to take the endpoints of the standard arcs 
$L_+\sqcup L_-$ (described in Step 1) to the endpoints of $L'$.  Now, by 
Lemma~\ref{facts}, $f$ is isotopic to a contactomorphism.  The isotopy is fixed 
on $M_1^c$ and might move one of the boundary components of $A$, but this can be 
extended over $B_s.$  Thus we can identify the tight contact manifolds $M_1$ and 
$M_1'$.  Moreover, according to Lemma~\ref{facts}, $S_+^m\circ S_-^n(L_1)$ is 
Legendrian isotopic to $L'_1$ and  $S_+^{-m}\circ S_-^{-n}(L_2)$ is Legendrian 
isotopic to $L_2'$, where $m$, $n$ are nonnegative integers. \end{proof}

\n
{\bf Step 3.} We now finish the proof of Theorem~\ref{mainambient} by using the 
following Lemma~\ref{lem3} to reduce to the previous step.

\begin{lem}\label{lem3}
Let $S$, $S'$ be (smoothly) isotopic standard convex separating spheres for 
$(M,L)$, with $S\cap N=S'\cap N$. Then there exists a finite sequence $S_0=S, 
S_1,\dots,S_l=S'$ of standard convex separating spheres where, for 
$i=0,\dots,k-1$, the pair $(S_i, S_{i+1})$ cobounds a region diffeomorphic to 
$S^2\times[0,1]$. 
\end{lem}

\begin{proof}[Proof of Lemma~\ref{lem3}.]
We use Colin's {\em isotopy discretization} technique \cite{Colin}.  Let 
$\Sigma_t$, $t\in[0,1]$, be the images of a smooth isotopy which takes 
$\Sigma_0=S$ to $\Sigma_1=S'$.  We may additionally assume that each $\Sigma_t$ 
intersects the standard neighborhood $N$ of $L$ in meridional ruling curves and 
that each $\Sigma_t\cap N$ is a standard convex meridional disk.  For each $t$, 
there exists a tubular neighborhood $N(\Sigma_t)$ of $\Sigma_t$ and an interval 
$[t-\varepsilon,t+\varepsilon]$  such that $\Sigma_s\subset N(\Sigma_t)$ for all 
$s\in [t-\varepsilon,t+\varepsilon]$.  By the compactness of $[0,1]$, there 
exist $t_0=0< t_1<\dots<t_k=1$ such that each $\Sigma_{t_{i+1}}$ is contained in 
a tubular neighborhood $N(\Sigma_{t_i})$ of $\Sigma_{t_i}$.  Since convex 
surfaces are $C^\infty$-dense in the space of closed embedded surfaces 
\cite{Giroux}, we may assume that the $\Sigma_{t_i}$  and $\bdry 
(N(\Sigma_{t_i}))= \Sigma'_{t_i} -\Sigma''_{t_i}$ are convex.  Now simply take 
the sequence $$ \Sigma_{t_0}, \Sigma'_{t_0}, \Sigma_{t_1},\Sigma'_{t_1},\dots.$$ 
It is easily verified that this sequence satisfies the cobounding condition. 
\end{proof}

\section{Proof of Theorem~\ref{mainsum}}\label{Proofs}

In this section we complete the proof of Theorem~\ref{mainsum}.  For simplicity 
we assume that $\mathcal{K}= \mathcal{K}_1\# \mathcal{K}_2$ and there are no 
equivalence relations of type 2 in Theorem~\ref{mainsum}, {\em i.e.}, there are 
no extra symmetries.  We show that 
\begin{equation} 
C:\left(\frac{\mathcal{L}(\mathcal{K}_1)\times\mathcal{L}(\mathcal{K}_2)}{\sim}
\right)\to \mathcal{L}(\mathcal{K}_1\#\mathcal{K}_2) 
\end{equation}
given by $(L_1,L_2)\mapsto L_1\# L_2$ is a bijection.  The proof is broken down 
into the following three claims. 

\begin{claim}
	The connect sum operation is well-defined.
\end{claim}

\begin{proof}
Theorem~\ref{maintech} indicates the ambiguity in splitting a manifold along 
different standard convex separating spheres.  Since we are always removing a 
standard trivial pair $(B_s,L_s)$ from a manifold/knot pair when taking a 
connected sum, and $L_s$ is not stabilized, the only state transition that 
actually occurs (among the possibilities listed in Theorem~\ref{maintech}) is 
Legendrian isotopy.  Hence there is no ambiguity in the connected sum operation. 
\end{proof}

It is clear that $C(S_\pm(L_1),L_2)$ and $C(L_1,S_\pm (L_2))$ are isotopic 
Legendrian knots, since stabilizations of a knot can be transfered from one side 
of the separating sphere to the other.  Therefore the map $C$ is well-defined.

\begin{claim}
	The map $C$ is surjective. 
\end{claim}

\begin{proof}
Let $L$ be a Legendrian knot in $\mathcal{L}(\mathcal{K}_1\# \mathcal{K}_2)$ and 
$S$ be a 2-sphere in $M^3$ that intersects $L$ transversely in exactly two 
points and topologically divides the knot into the appropriate knot types. 
Also let $N$ be a standard convex neighborhood of $L$, where we arrange 
the ruling curves on $\partial N$ to be meridional.  First isotop $S$ 
so that $S\cap \partial N$ consists of precisely two ruling curves, and then 
apply a further isotopy of $S$ rel $S\cap \partial N$ so that $S$ becomes 
convex.  Denote by $(M^o_1, L^o_1)$ and $(M^o_2,L^o_2)$ the 
components of the cut-open manifold $M\setminus S$ together with the cut-open 
Legendrian knot $L\setminus S$.  Let $L_s$ be a trivial Legendrian arc in 
$B_s$.  We now glue the standard contact 3-ball $B_s$ (with convex boundary) 
onto $M^o_i$, $i=1,2$, to form a closed tight contact manifold $M_i$; at the 
same time we glue $L_s$ and $L^o_i$ into a Legendrian knot $L_i\subset 
M_i$ with $L_i\in \mathcal{L}(\mathcal{K}_i)$.  Moreover, since we formed the 
connected sum of $L_1$ and $L_2$ by removing $B_s$ from each of $M_1$ 
and $M_2$ and gluing the resulting boundaries together, it is also clear that 
$L=L_1\# L_2.$ \end{proof}

\begin{claim}
	If $C(L_1, L_2)=C(L'_1, L'_2)$, then $(L_1, L_2)\sim (L'_1, L'_2).$ 
\end{claim}

\begin{proof}
Assume that $L_1\# L_2=L'_1\# L'_2,$ and let $S$ (resp.\ $S'$) be a standard 
convex separating sphere for $L_1\# L_2$ (resp.\ $L'_1\# L'_2$.)  Since $S$ and 
$S'$ are smoothly isotopic, Theorem~\ref{maintech} implies that $(L_1, L_2)\sim 
(L'_1, L'_2).$ \end{proof}

\s\s
\noindent
{\em Acknowledgments.}
Much of this work was done during the Workshop on Contact Geometry in the fall 
of 2000.   The authors gratefully acknowledge the support of Stanford 
University and the American Institute of Mathematics during this workshop.  The 
first author was supported by an NSF Postdoctoral Fellowship (DMS-9705949) and by
NSF Grand \# DMS-0203941, and 
the second by NSF Grant \# DMS-0072853 and an Alfred P.\ Sloan Research 
Fellowship. The first author also expresses his gratitude for the hospitality and 
generosity of Professor T.~Tsuboi and University of Tokyo
during the completion of this paper. 


\end{document}